\documentclass[twoside]{article}
\include{graphicx}
\title{Comments on ``Exactification of Stirling's approximation for the logarithm of the gamma function''}
\author{R. B. Paris\\
\\
School of Engineering, Computing and Applied Mathematics,\\
 University of Abertay Dundee, Dundee DD1 1HG, UK\\
E-Mail: r.paris@abertay.ac.uk}
\begin{document}
\def\f#1#2{\mbox{${\textstyle \frac{#1}{#2}}$}}
\def\dfrac#1#2{\displaystyle{\frac{#1}{#2}}}
\def\boldal{\mbox{\boldmath $\alpha$}}
\newcommand{\bee}{\begin{equation}}
\newcommand{\ee}{\end{equation}}
\newcommand{\lam}{\lambda}
\newcommand{\ka}{\kappa}
\newcommand{\al}{\alpha}
\newcommand{\om}{\omega}
\newcommand{\Om}{\Omega}
\newcommand{\fr}{\frac{1}{2}}
\newcommand{\fs}{\f{1}{2}}
\newcommand{\g}{\Gamma}
\newcommand{\br}{\biggr}
\newcommand{\bl}{\biggl}
\newcommand{\ra}{\rightarrow}
\renewcommand{\topfraction}{0.9}
\renewcommand{\bottomfraction}{0.9}
\renewcommand{\textfraction}{0.05}
\newcommand{\mcol}{\multicolumn}
\date{}
\maketitle
\begin{abstract}
We re-examine the exponentially improved expansion for $\log\,\g(z)$, first considered in Paris and Wood in 1991, to point out that the recent treatment by Kowalenko [Exactification of Stirling's approximation for the logarithm of the gamma function, arXiv:1404.2705] using his procedure of regularisation produces an equivalent result. 
In addition, we point out an error in his definition of the Stokes multiplier that leads him to make the incorrect statement that the Stokes phenomenon is a jump discontinuity, rather than a smooth transition. We supply a numerical example that clearly demonstrates the smooth transition of the leading subdominant exponential $e^{2\pi iz}$ across the Stokes line $\arg\,z=\fs\pi$.

\vspace{0.4cm}

\noindent {\bf Mathematics Subject Classification:} 34E05, 30C15, 30E15, 34E05, 41A60

\vspace{0.3cm}

\noindent {\bf Keywords:} Exponentially improved expansion, Stokes phenomenon, Gamma function
\end{abstract}
\vspace{0.3cm}

\begin{center}
{\bf 1. \  Introduction}
\end{center}
\setcounter{section}{1}
\setcounter{equation}{0}
\renewcommand{\theequation}{\arabic{section}.\arabic{equation}}
It is well known (see, for example, \cite[p.~277]{WW}) that as $|z|\ra\infty$ in the sector $|\arg\,z|\leq\pi-\epsilon$, $\epsilon>0$, the logarithm of the gamma function has the asymptotic expansion
\bee\label{e11}
\log\,\g(z)\sim (z-\fs)\log\,z-z+\fs\log\,2\pi+\sum_{r=1}^\infty\frac{B_{2r}}{2r(2r-1)z^{2r-1}},
\ee
where $B_{2r}$ are the Bernoulli numbers. Since successive even-order Bernoulli numbers have opposite signs, all terms in the asymptotic series have the same phase on $\arg\,z=\pm\fs\pi$, with the consequence that the positive and negative imaginary axes are Stokes lines.

In \cite{PW}, Paris and Wood obtained a refinement of this result that involved a finite truncation after $N$ terms of the above asymptotic series together with a remainder term $R_N(z)$ expressed as a convergent infinite sum of exponentials $e^{\pm2\pi ikz}$, $k=1, 2, \ldots $, with coefficients given by incomplete gamma functions.
If the slowly varying part of $\log\,\g(z)$ is defined by
\bee\label{e12}
\Om(z):=\log\,\g(z)-(z-\fs)\log\,z+z-\fs\log\,2\pi,
\ee
they found that \cite[Eq. (4.1), (4.11)]{PW}
\bee\label{e13}
\Om(z)=\sum_{r=1}^{N-1} \frac{B_{2r}}{2r(2r-1)z^{2r-1}}+R_N(z),
\ee
where
\bee\label{e14}
R_N(z)=\sum_{k=1}^\infty \frac{1}{k}\{e^{2\pi ikz}T_\nu(2\pi ikz)-e^{-2\pi ikz}T_\nu(-2\pi ikz)\}, \qquad \nu=2N-1
\ee
and $T_\nu(z)$ is the so-called terminant function defined as a multiple of the incomplete gamma function $\g(a,z)$ by
\bee\label{e15}
T_\nu(z)=\frac{e^{\pi i\nu}\g(\nu)}{2\pi i}\,\g(1-\nu,z).
\ee

The exponentials in $R_N(z)$ are subdominant relative to the finite series in (\ref{e13}) in the upper and lower half-planes and are maximally subdominant on the Stokes lines $\arg\,z=\pm\fs\pi$, respectively. It was established in \cite{PW} that, when the finite series is optimally truncated at or near its smallest term, the coefficients (the Stokes multipliers) of the leading subdominant exponentials (corresponding to $k=\pm 1$) undergo a smooth transition in the neighbourhood of the Stokes lines $\arg\,z=\pm\fs\pi$ given approximately by
\bee\label{e6}
\fs\pm\fs\mbox{erf}\,[(\theta\mp\fs\pi)\sqrt{\pi |z|}],\qquad \theta=\arg\,z.
\ee
This follows the error-function smoothing law first 
developed by Berry \cite{B1}. Subsequently, Berry \cite{B2} showed, by a sequence of increasingly delicate subtractions of optimally truncated asymptotic series, that all the exponentials switch on smoothly across the Stokes lines in a similar manner; see also the account given in \cite[\S 6.4]{PK}. 

In a recent and very long paper, Kowalenko \cite{Kow} employed his theory of regularisation 
to derive an expansion for $\log\,\g(z)$ that is equivalent to (\ref{e13}) and (\ref{e14}). He also expressed the remainder $R_N(z)$ as 
an infinite sum of integrals, which follows from (\ref{e14}) by expressing the incomplete gamma functions in their equivalent integral form. He carried out a detailed numerical study of his expansion to demonstrate that
$\log\,\g(z)$ can be computed to high accuracy (50dp) independently of the choice of the truncation index $N$. As we show in Section 3 these numerical results are basically correct. However, he goes on to make the assertion that, based on his numerics, the transition of the Stokes multiplier across the Stokes lines is not smooth, but discontinuous
(as originally proposed by Stokes) jumping from 0 to 1 as one crosses $\arg\,z=\pm\fs\pi$ in the sense of increasing $|\arg\,z|$. Indeed, he goes further to say \cite[p.~29]{Kow} that proponents of the smooth-transition theory of the Stokes multiplier have never displayed any numerical evidence to support their arguments. This is manifestly false as numerical results confirming the smooth error-function transition of the Stokes multiplier have been given, for example,  in \cite{B1} and \cite[pp.~257--9, 288]{PK}.

The purpose of this note is twofold: (i) to reiterate that the expansion (\ref{e13}) can be derived without recourse to the notion of regularisation and, more importantly, (ii) to point out the erroneous nature of Kowalenko's definition of the Stokes multiplier that leads him to his incorrect conclusions. In the final section, we conclude with a numerical example that clearly demonstrates the error-function smoothing of the Stokes multiplier associated with the leading subdominant exponential $e^{2\pi iz}$ in $R_N(z)$ across $\arg\,z=\fs\pi$.

\vspace{0.6cm}

\begin{center}
{\bf 2. \ The derivation of the exponentially improved expansion}
\end{center}
\setcounter{section}{2}
\setcounter{equation}{0}
\renewcommand{\theequation}{\arabic{section}.\arabic{equation}}
We consider $\theta=\arg\,z\geq 0$, since values of $\log\,\g(z)$ corresponding to $\theta<0$ take conjugate values. Also, unlike Kowalenko \cite{Kow} who treats $\arg\,z$ lying in sectors beyond $\pm\pi$, we restrict our attention to the sector $0\leq\theta<\pi$ since this includes the region of prime interest, namely the Stokes line $\theta=\fs\pi$.

The analysis we present in this section is essentially that given in Paris and Wood \cite{PW} which we repeat here for completeness in exposition. We start with a Mellin integral representation for the slowly varying part $\Om(z)$ of the logarithm of the gamma function in (\ref{e13}) given by \cite[pp.~277--278]{WW}, \cite[p.~282]{PK}
\[\Om(z)=-\frac{1}{2\pi i}\int_{c-\infty i}^{c+\infty i}\frac{\pi z^{-s}\zeta(-s)}{s \sin \pi s}\,ds\qquad (0<c<1)\]
valid when $|\arg\,z|\leq\pi-\epsilon$, where $\zeta(s)$ is the Riemann zeta function. If we employ the functional relation for $\zeta(s)$ in the form \cite[p.~269]{WW}
\[\pi\zeta(-s)=-(2\pi)^{-s} \g(1+s)\zeta(1+s) \sin \fs\pi s,\]
followed by displacement of the integration path to the right over the simple poles of the integrand at $s=1, 3, \ldots , 2N-1$, where $N$ is an arbitrary positive integer, we obtain
\begin{eqnarray}
\Om(z)&=&\frac{1}{2\pi i}\int_{c-\infty i}^{c+\infty i} (2\pi z)^{-s}\g(1+s)\zeta(1+s) \frac{\sin \fs\pi s}{s \sin \pi s}\,ds\nonumber\\
&=&\sum_{r=1}^{N-1}\frac{B_{2r}}{2r(2r-1)z^{2r-1}}+R_N(z),\label{e21}
\end{eqnarray}
where we have used the result connecting the even-index Bernoulli numbers to the $\zeta$ function $B_{2r}=2(-)^{r-1}(2r)! \zeta(2r)/(2\pi)^{2r}$. The remainder term $R_N(z)$ is given by
\begin{eqnarray*}
R_N(z)&=&\frac{1}{2\pi i}\int_{L_N}(2\pi z)^{-s}\,\frac{\g(s)}{\sin \pi s}\,\zeta(1+s) \sin \fs\pi s\,ds\\
&=&\frac{1}{2\pi i}\sum_{k=1}^\infty \frac{1}{k}\int_{L_N}(2\pi kz)^{-s}\,\frac{\g(s)}{\sin \pi s}\,\sin \fs\pi s\,ds,
\end{eqnarray*}
where $L_N$ denotes the displaced integration path $(-\infty i+2N-1-c, \infty i+2N-1-c)$ with $0<c<1$ and
we have employed the Dirichlet series expansion for $\zeta(1+s)$ (which is permissible since Re $(s)>1$ on the displaced integration path when $N\geq 1$).

Writing $\sin \fs\pi s$ in terms of exponentials, we can express the remainder after $N$ terms in the form
\bee\label{e22}
R_N(z)=\frac{1}{2i}\sum_{k=1}^\infty \frac{1}{k} (J_--J_+),
\ee
where,\footnote{Although the integrals $J_\pm$ are valid in $|\arg\,(\pm iz)|<\frac{3}{2}\pi$, the combination $J_--J_+$ has the common sector of validity $|\arg\,z|<\pi$.} with $\chi\equiv 2\pi kz$,
\begin{eqnarray}
J_\pm=\frac{1}{2\pi i}\int_{L_N} (\pm i\chi)^{-s}\,\frac{\g(s)}{\sin \pi s}\,ds
&=&-\frac{(\pm i\chi)^{1-2N}}{2\pi i}\int_{-c-\infty i}^{-c+\infty i} (\pm i\chi)^{-\tau}\,\frac{\g(\tau+2N-1)}{\sin \pi\tau}\,d\tau\nonumber\\
&=&\frac{\g(2N-1)}{\pi}\,e^{\pm i\chi} \g(2-2N,\pm i\chi).\label{e23}
\end{eqnarray}
In obtaining the last expression we have made the change of variable $s=\tau+2N-1$ and
employed the standard result involving the incomplete gamma function $\g(a,z)$ \cite[Eq.~(8.6.12)]{DLMF}, \cite[p.~113]{PK}
\bee\label{e24}
\frac{1}{2\pi i}\int_{-c-\infty i}^{-c+\infty i} \frac{\g(s+\nu)}{\sin \pi s} \,z^{-s}ds=-\frac{\g(\nu)}{\pi} z^\nu e^z\,\g(1-\nu,z) \qquad (|\arg\,z|<\f{3}{2}\pi),
\ee
where the integration path (indented if necessary) has $0<c<1$ and passes to the right of the poles of $\g(s+\nu)$.

Then, from (\ref{e22}) and (\ref{e23}), we finally obtain
\bee\label{e25}
R_N(z)=-\frac{\g(2N-1)}{2\pi i} \sum_{k=1}^\infty \frac{1}{k}\{e^{2\pi ikz} \g(2-2N,2\pi ikz)-e^{-2\pi ikz} \g(2-2N,-2\pi ikz)\}.
\ee
From the asymptotic behaviour $\g(a,x)\sim x^{a-1}e^{-x}$ as $|x|\ra\infty$ in $|\arg\,x|<\f{3}{2}\pi$ \cite[p.~179]{DLMF}, it can be seen that the late terms in the infinite series on the right-hand side of (\ref{e25}) are O$(k^{-2N})$ as $k\ra\infty$ when $|\arg\,z|\leq\pi-\epsilon$. The sum in (\ref{e25}) is therefore absolutely convergent when $N\geq 1$ in this sector. This last result combined with (\ref{e21}) then gives the exponentially improved expansion of $\Om(z)$ first obtained in Paris and Wood \cite[Eqs. (4.1), (4.11)]{PW}.
Another derivation starting from Binet's second representation for the logarithm of the gamma function is given in the appendix.

An alternative form for $R_N(z)$ is given in (\ref{a1}) as
\bee\label{e26}
R_N(z)=\frac{2(-)^{N-1}z}{(2\pi z)^{2N-2}}\sum_{k=1}^\infty\frac{1}{k^{2N-2}}\int_0^\infty \frac{w^{2N-2}e^{-w}}{w^2+4\pi^2k^2z^2}\,dw
\ee
after a simple change of variable, which holds in $|\arg\,z|\leq\fs\pi-\epsilon$. This last formula also follows immediately from (\ref{e25}) upon use of the integral representation (\ref{a3}) for the incomplete gamma function.

When computing $\g(a,x)$, {\it Mathematica} returns the value only in the principal sector $-\pi<\arg\,x\leq\pi$. 
Consequently, we can compute $R_N(z)$ in $|\arg\,z|\leq\fs\pi$ using (\ref{e25}), {\it thereby including the Stokes lines} $\arg\,z=\pm\fs\pi$. However, when we compute $R_N(z)$ in the sector $\fs\pi<\arg\,z<\pi$, we either have to use the analytic continuation for $\g(2-2N,2\pi ikz)$ given in \cite[Eq.~(8.2.10)]{DLMF} or, equivalently, make use of the continuation formula \cite[p.~281]{PK}
\[\Om(z)+\Om(ze^{\mp\pi i})=-\log\,(1-e^{\pm2\pi iz}),\]
which follows from the reflection formula for the gamma function $\g(z)=-\pi/(z \sin \pi z\, \g(-z))$.

Then, we can write
\bee\label{e27}
\Om(z)=\left\{\begin{array}{ll}\displaystyle{\sum_{r=1}^{N-1} \frac{B_{2r}}{2r(2r-1)z^{2r-1}}+R_N(z) } & (0\leq\theta\leq\fs\pi)\\
\\
\displaystyle{\sum_{r=1}^{N-1}\frac{B_{2r}}{2r(2r-1)z^{2r-1}}+R_N'(z)-\log\,(1-e^{2\pi iz}) } & (\fs\pi<\theta<\pi),\end{array}\right.
\ee
where the prime on $R_N(z)$ in the second expression indicates that the argument of the first incomplete gamma function in (\ref{e25}) is to be replaced\footnote{When computing $\g(a,2\pi ikz)$ with $\fs\pi<\arg\,z<\pi$, {\it Mathematica} automatically returns the desired value $\g(a,2\pi ikze^{-2\pi i})$.}
 by $2\pi ikze^{-2\pi i}$.
This is essentially the result obtained by Kowalenko \cite[Eq.~(72)]{Kow} except that he gives a separate expression on the Stokes line $\theta=\fs\pi$, which we have found not necessary when using (\ref{e25}) 
for the reason stated above. Kowalenko gave both forms (\ref{e25}) and (\ref{e26}) for the remainder but appears to have employed the expression (\ref{e26}) in his computations. We repeat that (\ref{e27}) has been derived by routine analysis which has not required regularisation as used in \cite{Kow}.  
 
The result (\ref{e27}) is exact. Consequently, within reasonable limits, arbitrary precision may be achieved
for $\log\,\g(z)$ irrespective of the value of $z$ and the truncation index $N$. These assertions are borne out by Kowalenko's calculations. It
does not matter whether the index $N$ is chosen to be the optimal truncation value $N_o\sim\pi |z|$ (\cite[p.~141]{PW}) or significantly different from this value, although its choice affects the rate of convergence of the series (\ref{e25}) for $R_N(z)$.
Kowalenko also applied (\ref{e27}) to the case of small $|z|$ (he took an extreme case where $z=10^{-1}$), where there is no optimal truncation index ($N=1$). As he points out, such a low value of $N$ makes it more difficult to compute $R_N(z)$ to a given accuracy on account of the decay of the late terms being controlled by $k^{-2N}$. If one chooses a larger value of $N$, the finite series in (\ref{e27}) becomes large for small $|z|$ and one is then confronted with the cancellation of large terms that results in loss of precision. However, it should be said that the main interest in (\ref{e27}) is for large values of $|z|$.

In Table 1 we present the absolute error in the computation of $\Om(z)$ from (\ref{e27}) compared with the value obtained from {\it Mathematica} using the LogGamma function. We take $z=5e^{i\theta}$, for which the optimal truncation index is $N_0=16$, and adjust the truncation index $k=K$ of the series for $R_N(z)$ in (\ref{e25}) to 
obtain a prescribed accuracy. The table shows values of $\theta$ in the right- and left-hand planes and on the Stokes line $\theta=\fs\pi$; in each case we obtain values accurate to over 50dp with the selected $K$ values.
The value $N=40$ represents a highly non-optimal case and corresponds to a situation where the final terms in the finite sum over $r$ rise to almost unity in magnitude.
\begin{table}[t]
\begin{center}
\begin{tabular}{l|ccc}
\mcol{1}{c|}{Truncation indices} & \mcol{1}{c}{$\theta=\f{1}{3}\pi$} & \mcol{1}{c}{$\theta=\fs\pi$} &  \mcol{1}{c}{$\theta=\f{3}{4}\pi$}\\
[.1cm]\hline
&&&\\[-0.25cm]
$N=12,\ K=40$ & $6.025\times 10^{-52}$ & $6.026\times10^{-52}$ & $6.024\times 10^{-52}$\\
$N=16,\ K=13$ & $7.789\times 10^{-52}$ & $7.809\times 10^{-52}$ & $7.770\times 10^{-52}$\\
$N=20,\ K=7$  & $5.226\times 10^{-51}$ & $5.293\times 10^{-51}$ & $5.161\times 10^{-51}$\\
$N=40,\ K=3$  & $1.487\times 10^{-52}$ & $2.244\times 10^{-52}$ & $1.210\times 10^{-52}$\\
[.1cm]\hline
\end{tabular}
\caption{\footnotesize{The absolute error in the computation of $\Om(z)$ using the expansion (\ref{e27}) for different $\theta$ when $z=5e^{i\theta}$. The optimal truncation index is $N_o=16$.}}
\end{center}
\end{table}

Finally, it is worth remarking that Kowalenko, who used the expansion (\ref{e26}) for $R_N(z)$ and so performed a series of numerical integrations, reported computation times running to several hours for each value of $z$. Our computation of $R_N(z)$ using the series of incomplete gamma functions in (\ref{e25}) took a fraction of a second to compute. 
 
\vspace{0.6cm}

\begin{center}
{\bf 3. \ The Stokes phenomenon }
\end{center}
\setcounter{section}{3}
\setcounter{equation}{0}
\renewcommand{\theequation}{\arabic{section}.\arabic{equation}}
Kowalenko \cite[p.~28]{Kow} called the term appearing in the second expression in (\ref{e27})
\[-\log\,(1-e^{2\pi iz})=\sum_{k=1}^\infty \frac{e^{2\pi ikz}}{k}\qquad (0<\arg\,z<\pi)\]
the `Stokes discontinuity term'. If this term is multiplied by the quantity $A\equiv A(\theta)$, then the multiplier $A$ undergoes (at fixed $|z|$) a step discontinuity across the Stokes line $\arg\,z=\fs\pi$, possessing the values\footnote{Kowalenko gave the values $A=0$ ($0\leq\arg\,z<\fs\pi$), $\fs$ ($\arg\,z=\fs\pi$), 1 ($\fs\pi<\arg\,z<\pi$).}
 0 and 1 in $0\leq\arg\,z\leq\fs\pi$ and $\fs\pi<\arg\,z<\pi$, respectively. His numerical computations using (\ref{e27}), of course, confirm this jump discontinuity. However, he {\it mistakenly calls the multiplier $A$ the Stokes multiplier associated with $\log\,\g(z)$} \cite[p.~29]{Kow}. This leads him to make the fallacious assertion that there is no smoothing of the Stokes phenomenon and that the change in this multiplier is a discontinuous jump as originally believed since Stokes' time.

The Stokes multiplier, which we shall denote by $S(\theta)$, is associated with a {\it single} subdominant exponential that is born (or extinguished) when crossing a Stokes line. The case of $\log \,\g(z)$ is unusual in that it involves not one but an infinite number of subdominant exponentials (see (\ref{e25})), each associated with its own Stokes multiplier that switches on across $\arg\,z=\fs\pi$. The treatment of these multipliers for the exponentials $e^{\pm2\pi ikz}$ ($k\geq 1$) is considered in \cite{B2}; see also the discussion in \cite[\S 6.4.2]{PK}. As stated in Section 1, we confine our attention here to the leading exponential $e^{2\pi iz}$
appearing in (\ref{e25}) and examine its behaviour for large $|z|$ in the neighbourhood of the Stokes line $\arg\,z=\fs\pi$.

The Stokes multiplier $S(\theta)$ for the leading subdominant exponential is {\it defined\/} (at fixed $|z|$) by
\bee\label{e31}
\Om(z)=\sum_{r=1}^{N_o-1}\frac{B_{2r}}{2r(2r-1)z^{2r-1}}+S(\theta)\,e^{2\pi iz},
\ee
where the finite series is optimally truncated at $N_0\simeq\pi|z|$. To obtain its value numerically, we therefore compute
\bee\label{e32}
S(\theta)=e^{-2\pi iz}\left\{\Om(z)-\sum_{r=1}^{N_o-1}\frac{B_{2r}}{2r(2r-1)z^{2r-1}}\right\}
\ee
for a series of $\theta$-values.
The approximate form of $S(\theta)$ in the neighbourhood of $\theta=\fs\pi$ can be deduced by employing the terminant function defined in (\ref{e15}) in the series (\ref{e25}) to yield the equivalent form of $R_N(z)$ given in (\ref{e14}). Retaining only the $k=1$ terms, we then have
\bee\label{e33}
R_{N_o}(z) \simeq e^{2\pi iz}T_\nu(2\pi iz)-e^{-2\pi iz}T_\nu(-2\pi iz), \qquad \nu=2N_o-1.
\ee

The asymptotics of $T_\mu(x)$ when $\mu\simeq|x|\ra\infty$ has been discussed in detail by Olver \cite{O}; see also
\cite[p.~67]{DLMF} and \cite[\S 6.2.6]{PK}. We have the leading terms, with $\phi=\arg\,x$,
\[T_\mu(x)\sim\left\{\begin{array}{lr}-\dfrac{ie^{(\pi-\phi)i\mu}}{1+e^{-i\phi}}\,\frac{e^{-x-|x|}}{\sqrt{2\pi |x|}} & -\pi+\epsilon\leq\phi\leq\pi-\epsilon\\
\\
\fs+\fs \mbox{erf}\,[c(\phi)(\fs |x|)^\fr]-\dfrac{iB_0e^{-\fr |x|c^2(\phi)}}{\sqrt{2\pi |x|}} & \epsilon\leq\phi\leq2\pi-\epsilon,
\end{array}\right.\]
where
\[\fs c^2(\phi)=1+i(\phi-\pi)-e^{i(\phi-\pi)},\qquad B_0=\frac{e^{(\pi-\phi)i\alpha}}{1+e^{-i\phi}}+\frac{i}{c(\phi)},\qquad \mu=|x|+\alpha,\]
with the branch for $c(\phi)$ chosen so that $c(\phi)\simeq \phi-\pi$ when $\phi\simeq \pi$, and 
the quantity $\alpha$ is bounded. When $\phi=\pi$ the limiting value of $B_0$ is $\f{2}{3}-\alpha$.

Then, with $x=\pm 2\pi iz$, $\phi=\theta\pm\fs\pi$ and $\mu=2N_o-1$ we obtain from (\ref{e31}) and (\ref{e33})
\bee\label{e34}
R_{N_o}(z)\sim e^{2\pi iz}\,S(\theta),\qquad S(\theta)\sim\fs+\fs \mbox{erf}\,[c(\theta+\fs\pi)\sqrt{\pi |z|}]-\frac{iC_0e^{-2\pi\gamma |z|}}{2\pi \sqrt{|z|}},
\ee
where, with $\om:=\theta-\fs\pi$,
\[\gamma=1+ie^{i\theta},\qquad C_0=B_0e^{-2\pi i\om |z|}+\frac{e^{-i\om\nu}}{1+e^{-i\om}},
\qquad B_0=\frac{e^{-i\om\alpha}}{1-e^{-i\om}}+\frac{i}{c(\theta+\fs\pi)},\]
\[\alpha=2N_0-1-2\pi|z|,\qquad c(\theta+\fs\pi)=\om+\f{1}{6}i\om^2-\f{1}{36}\om^3+\f{1}{270}i\om^4+\cdots\ .\]
We note that $\gamma=0$ and $C_0=\f{7}{6}-\alpha$ on the Stokes line $\theta=\fs\pi$. Although the factor $e^{-2\pi\gamma|z|}$ in the approximation for $S(\theta)$ in (\ref{e34}) decays exponentially away from this ray, it makes a non-negligible contribution, particularly to Im\,$S(\theta)$, in the immediate vicinity of the Stokes line.

Since $c(\theta+\fs\pi)\simeq \theta-\fs\pi$ near 
$\arg\,z=\fs\pi$, we find from (\ref{e34}) that approximately
\bee\label{e35}
\mbox{Re}\,S(\theta)\simeq \fs+\fs\mbox{erf}\,[(\theta-\fs\pi) \sqrt{\pi|z|}] 
\ee
as first found in Paris and Wood \cite{PW}, with Im\,$S(\theta)$ given by the imaginary part of (\ref{e34}).
Values of the real and imaginary parts of $S(\theta)$ obtained from (\ref{e32}) and the approximation (\ref{e34}) when $z=8e^{i\theta}$ as a function of $\theta$ are presented in Table 2. These values are also illustrated in Fig.~1. 
\begin{figure}[hb]
	\begin{center}
		{\tiny($a$)}\includegraphics[width=0.40\textwidth]{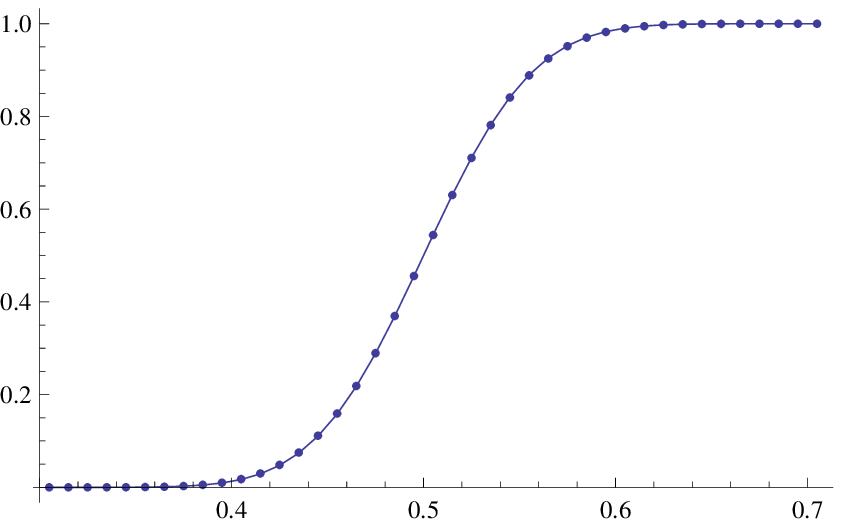}\hspace{1.2cm} {\tiny($b$)}\includegraphics[width=0.40\textwidth]{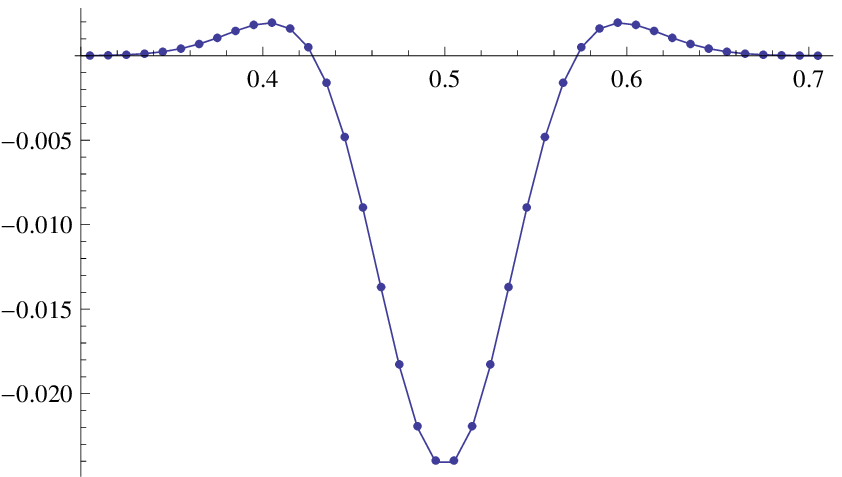}
	\caption{\small{The variation of (a) Re $S(\theta)$ and (b) Im $S(\theta)$ as a function of $\theta/\pi$ when $z=8e^{i\theta}$ and the optimal truncation index $N_o=26$. The curves represent the approximate values from (\ref{e34}) and the dots the exact values from (\ref{e32}). }}\label{f1}
	\end{center}
\end{figure}
It is seen that $S(\theta)$ has a small imaginary part that is extremely well accounted for by (\ref{e34}). Values of $S(\theta)$ in the case $|z|=5$ are given in \cite[p.~288]{PK} and are compared with those given by the formula (\ref{e35}).
\begin{table}[t]
\begin{center}
\begin{tabular}{l|c|c}
\mcol{1}{c|}{$\theta/\pi$} & \mcol{1}{c|}{$S(\theta)$} & \mcol{1}{c}{Approx. $S(\theta)$} \\
[.1cm]\hline
&\\[-0.25cm]
0.325 & $0.0000224+0.0000580i$ & $0.0000223+0.0000580i$\\
0.350 & $0.0003638+0.0003179i$ & $0.0003641+0.0003176i$\\
0.400 & $0.0133877+0.0019356i$ & $0.0133896+0.0019274i$\\
0.450 & $0.1338254-0.0067919i$ & $0.1338262-0.0068425i$\\
0.475 & $0.2894310-0.0182669i$ & $0.2894307-0.0183467i$\\
0.500 & $0.5000000-0.0242241i$ & $0.5000000-0.0243169i$\\
0.525 & $0.7105689-0.0182669i$ & $0.7105693-0.0183467i$\\
0.550 & $0.8661746-0.0067919i$ & $0.8661738-0.0068425i$\\
0.600 & $0.9866123+0.0019356i$ & $0.9866104+0.0019274i$\\
0.650 & $0.9996362+0.0003179i$ & $0.9996359+0.0003176i$\\
0.700 & $1.0000017+0.0000060i$ & $1.0000000+0.0000060i$\\
0.750 & $1.0000000+0.0000000i$ & $1.0000000+0.0000000i$\\
[.1cm]\hline
\end{tabular}
\caption{\footnotesize{Values of the Stokes multiplier $S(\theta)$ for different $\theta/\pi$ when $z=8e^{i\theta}$
compared with the approximation (\ref{e34}). The optimal truncation index is $N_o=26$ and the parameter $\alpha\doteq 0.734518$.}}
\end{center}
\end{table}
\vspace{0.6cm}

\begin{center}
{\bf 4. \ Concluding remarks }
\end{center}
\setcounter{section}{4}
\setcounter{equation}{0}
\renewcommand{\theequation}{\arabic{section}.\arabic{equation}}
We have repeated a calculation, first given in Paris and Wood \cite{PW}, which gives the exponentially improved expansion for $\log\,\g(z)$ valid in $|\arg\,z|\leq\pi-\epsilon$. This expansion has been obtained by standard analysis, both using a Mellin integral and Binet's representation for $\log\,\g(z)$. It agrees with the recent result obtained by Kowalenko \cite[Eq.~(72)]{Kow} using his process of regularisation.

We have pointed out that, although Kowalenko's numerical computations are basically correct, his interpretation of the Stokes multiplier for the Stokes lines $\arg\,z=\pm\fs\pi$ is incorrect. This has led him to make the statement that there is no smoothing of the Stokes phenomenon. A correct definition of the Stokes multiplier $S(\theta)$ associated with the leading subdominant exponential $e^{2\pi iz}$ appearing in the remainder term leads to the predicted error-function smoothing across the Stokes line $\arg\,z=\fs\pi$; an analogous result applies for the Stokes line $\arg\,z=-\fs\pi$. It has also been shown in a particular case that the approximate variation of the real and imaginary parts of $S(\theta)$ follow very closely the numerically computed values.

\vspace{0.6cm}

\begin{center}
{\bf Appendix: \ Alternative derivation of expansion for $\Om(z)$}
\end{center}
\setcounter{section}{1}
\setcounter{equation}{0}
\renewcommand{\theequation}{\Alph{section}.\arabic{equation}}
We use Binet's second representation for $\Om(z)$ given by \cite[p.~251]{WW}
\[\Om(z)=2\int_0^\infty \frac{\arctan (t/z)}{e^{2\pi t}-1}\,dt\]
valid when Re$(z)>0$, where $\arctan$ has its principal value. Using the expansion for arbitrary positive integer $N$
\[\arctan (t/z)=\sum_{r=1}^{N-1} \frac{(-)^{r-1} (t/z)^{2r-1}}{2r-1}+\frac{(-)^{N-1}}{z^{2N-1}}\int_0^t \frac{u^{2N-2}}{u^2+z^2}\,du,\]
together with the integral 
$\int_0^\infty t^{2r-1}/(e^{2\pi t}-1)\,dt=(-)^{r-1}B_{2r}/(4r)$,
we obtain
\[\Om(z)=\sum_{r=1}^{N-1}\frac{B_{2r}}{2r(2r-1) z^{2r-1}}+R_N(z),\]
where \cite[p.~252]{WW}
\[R_N(z)=\frac{2(-)^{N-1}}{z^{2N-3}}\int_0^\infty \!\!\!\int_0^t\frac{u^{2N-2}}{u^2+z^2}\,\frac{1}{e^{2\pi t}-1}\,du\,dt.\]

Expansion of the exponential factor followed by interchange in the order of integration then yields
\begin{eqnarray}
R_N(z)&=& \frac{2(-)^{N-1}}{z^{2N-3}}\sum_{k=1}^\infty \int_0^\infty \!\!\!\int_u^\infty \frac{u^{2N-2}e^{-2\pi kt}}{u^2+z^2}\,dt\,du\nonumber\\
&=&\frac{2(-)^{N-1}}{z^{2N-3}}\sum_{k=1}^\infty \frac{1}{2\pi k}\int_0^\infty \frac{u^{2N-2}e^{-2\pi ku}}{u^2+z^2}\,du.\label{a1}
\end{eqnarray}
If the denominator is decomposed into partial fractions, we can then use the integral representation for the incomplete gamma function \cite[Eq.~(8.6.4)]{DLMF}
\bee\label{a3}
\g(a,z)=\frac{z^ae^{-z}}{\g(1-a)}\int_0^\infty\frac{t^{-a}e^{-t}}{z+t}\,dt\qquad (|\arg\,z|<\pi,\ 
\mbox{Re} (a)<1)
\ee
to find after a simple change of variable
\bee\label{a2}
R_N(z)=-\frac{\g(2N-1)}{2\pi i} \sum_{k=1}^\infty \frac{1}{k}\{e^{2\pi ikz} \g(2-2N,2\pi ikz)-e^{-2\pi ikz} \g(2-2N,-2\pi ikz)\},
\ee
as obtained using the Mellin integral approach in Section 2. The result (\ref{a2}) has been established for Re $(z)>0$,
but may continued analytically into the wider sector $|\arg\,z|<\pi$, corresponding to the domain of analyticity of $\Om(z)$.

\end{document}